\documentclass{amsart}

\newtheorem{theorem}{Theorem}
\newcommand{\bt}{\begin{theorem}}
\newcommand{\et}{\end{theorem}}
\newtheorem{lemma}{Lemma}
\newcommand{\bl}{\begin{lemma}}
\newcommand{\el}{\end{lemma}}
\newtheorem{corollary}{Corollary}
\newcommand{\bc}{\begin{corollary}}
\newcommand{\ec}{\end{corollary}}
\newtheorem{problem}{Problem}
\newcommand{\bprob}{\begin{problem}}
\newcommand{\eprob}{\end{problem}}
\newtheorem*{conjectureNN}{Conjecture}
\newcommand{\bconjNN}{\begin{conjectureNN}}
\newcommand{\econjNN}{\end{conjectureNN}}
\newtheorem*{theoremNN}{Theorem}
\newcommand{\btNN}{\begin{theoremNN}}
\newcommand{\etNN}{\end{theoremNN}}
\newcommand{\beq}{\begin{equation}}
\newcommand{\eeq}{\end{equation}}
\newcommand{\benum}{\begin{enumerate}}
\newcommand{\eenum}{\end{enumerate}}
\newcommand{\N}{\ensuremath{ \mathbf N }}
\newcommand{\Z}{\ensuremath{\mathbf Z}}

\newcommand{\mcf}{\ensuremath{ \mathcal F}}
\newcommand{\mcg}{\ensuremath{ \mathcal G}}

\newcommand{\bmat}{\left(\begin{matrix}}
\newcommand{\emat}{\end{matrix}\right)}

\DeclareMathOperator{\qqand}{\qquad\text{and}\qquad}
\DeclareMathOperator{\card}{\text{card}}
\DeclareMathOperator{\kernel}{\text{kernel}}

\title[Commutative algebra and the Frobenius problem]
{Commutative algebra and the linear diophantine problem of Frobenius}

\author{Melvyn B. Nathanson}\address{Department of Mathematics\\Lehman College (CUNY)\\Bronx, NY 10468} \email{melvyn.nathanson@lehman.cuny.edu}


\begin{document}

\maketitle

\begin{abstract}  
Let $A$ be a finite set of relatively prime positive integers, and let $S(A)$ be the set 
of all nonnegative integral linear combinations of elements of $A$.  
The set $S(A)$ is a semigroup that contains all sufficiently large integers.  
 The largest integer not in $S(A)$ is the \emph{Frobenius number} of $A$, 
 and the number of positive integers not in $S(A)$ is the \emph{genus} of $A$.  
 Sharp and Sylvester proved in 1884 that the Frobenius number of the set $A = \{a,b\}$ is 
 $ab-a-b$, and that the genus of $A$ is $(a-1)(b-1)/2$.  
 Graded rings and a simple form of Hilbert's syzygy theorem are used to give a commutative 
 algebra proof of this result.  
 \end{abstract}

\section{The linear diophantine problem of Frobenius}

Let $\N_0$ be the additive semigroup of nonnegative integers.  
A \emph{numerical semigroup} is a subsemigroup $S$ of $\N_0$ 
that contains 0 and contains all sufficiently large integers.

\bt                     \label{FP:theorem:basic}
Let $A  = \{a_1,\ldots, a_k\}$  be a finite set of positive integers with $\card(A) = k \geq 2$, 
and let $\gcd(A)$ denote the greatest common divisor of the integers in $A$.
The set 
\beq      \label{FP:S(A)}
S(A) = \left\{ \sum_{i=1}^k a_i r_i : r_i \in \N_0 \text{ for all } i \in \{1,2,\ldots, k\} \right\}
\eeq
is a numerical semigroup if and only if $\gcd(A) = 1$. 
\et 

\begin{proof}
The set $S(A)$ is a semigroup because it contains 0 and is closed under addition.
Every integer in $S(A)$ is divisible by $\gcd(A)$.  
If $S(A)$ is a numerical semigroup, then it contains all sufficiently large integers, 
and so $\gcd(A) = 1$.

A basic theorem in elementary number theory  states that if $\gcd(A) = 1$, 
then every integer $n$ can be written as an integral linear combination 
of the elements of $A$.  If 
\[
n =  \sum_{i=1}^k a_i s_i 
\]
with $s_1, \ldots, s_k \in \Z$, then    
\[
n \equiv \sum_{i=1}^{k-1} a_i s_i \pmod{a_k}.
\]
For each $i \in \{1,\ldots, k-1\}$, there exists $r_i \in \{0,1,2,\ldots, a_k-1\}$ such that 
$
s_i \equiv r_i \pmod{a_k}   
$
and so
\[
n \equiv \sum_{i=1}^{k-1} a_i r_i \pmod{a_k}.
\]
There is an integer $r_k$ such that 
\[
a_kr_k = n - \sum_{i=1}^{k-1} a_i r_i  \geq  n - (a_k - 1) \sum_{i=1}^{k-1} a_i. 
\] 
If $n \geq (a_k - 1) \sum_{i=1}^{k-1} a_i$, then $r_k \geq 0$, 
and $n =  \sum_{i=1}^k a_i r_i $ is a representation of $n$ 
as a nonnegative integral linear combination of the elements of the $A$.
This completes the proof.  
\end{proof}

Let $A$  be a finite set of positive integers with $\gcd(A) = 1$.  
The \emph{Frobenius number} of $A$ is the largest integer $\mcf(A)$ 
not contained in $S(A)$.    
The proof of Theorem~\ref{FP:theorem:basic} shows that 
\[
\mcf(A) \leq  (a_k - 1) \sum_{i=1}^{k-1} a_i - 1.
\]
The elements of the finite set $\N_0 \setminus S(A) = \{0,1,2,\ldots, \mcf(A) \} \setminus S(A)$ 
are called the \emph{gaps}\index{gaps}\index{Frobenius!gaps} of $S(A)$.  
The \emph{genus}\index{genus}\index{Frobenius!genus}  of $A$, 
denoted $\mcg(A)$, is the number of gaps of $S(A)$.  
Because $S(A)$ is closed under addition and $\mcf(A) \notin S(A)$,  
it follows that if $n \in S(A)$, then  $\mcf(A) - n \notin S(A)$.  
Therefore, $S(A)$ contains at most one element of the set 
$\{n, \mcf(A)-n\}$ for all $n \in \{0,1,2,\ldots, \mcf(A)\}$, 
and so $\mcg(A) \geq (\mcf(A) + 1)/2$.
The  numerical semigroup $S(A)$ is \emph{symmetric} if  
$n \notin S(A)$ implies $\mcf(A) - n \in S(A)$.

The \emph{linear diophantine problem of Frobenius} is to compute the integer $\mcf(A)$.
In 1884, Sylvester~\cite{sylv84a} and Sharp~\cite{shar84} 
proved that the set  $A = \{a,b\}$ has Frobenius number 
$\mcf(a,b) = ab - a -b$ and  genus $\mcg(a,b) = (a-1)(b-1)/2$.

For sets $A$ with $|A| \geq 3$, the problem is still unsolved and mysterious.  
Indeed, there is no explicit solution to the Frobenius problem even for sets 
$A$ with $|A| = 3$.  
Methods from number theory, analysis, geometry, probability, and algebraic geometry  
have produced many partial results.  
Some of this is described in a monograph by Ramirez-Alfonsin~\cite{rami05}, and 
there is much recent work  (for example, Aliev-Henk~\cite{alie-henk09}, 
Arnold~\cite{arno06,arno07}, 
Bourgain-Sinai~\cite{bour-sina07}, 
Fel~\cite{fel08,fel09}, Fukshansky-Robins~\cite{fuks-robi07}, Marklof~\cite{mark10},  
Schmidt~\cite{schm15}, and Str{\"o}mbergsson~\cite{stro12}).  

In this paper  we show how elementary commutative algebra has been applied 
to obtain the Sharp-Sylvester solution to the Frobenius problem for $k = 2$.

\section{The Frobenius number $\mcf(a,b)$}

In Section~\ref{FP:section:proof}, we  use commutative algebra 
(essentially, a simple form of the Hilbert syzygy theorem) to prove the following result.  

\bt                           \label{FP:theorem:polyIdentity}
Let $A = \{a,b\}$, where $a$ and $b$ are distinct, relatively prime positive integers, 
and let $S(A) =  \{ai +bj: i,j \in \N_0\}$.  
The generating function for the gaps of the numerical semigroup $S(A)$ 
is the polynomial 
\beq                   \label{FP:f_A}
f_A(q) = \sum_{n \in \N_0 \setminus S(A)} q^n. 
\eeq 
This polynomial satisfies the functional equation
\beq        \label{FP:polyIdentity}
(q^a - 1) (q^b - 1) (  (q-1)f_A(q)+1 ) = (q-1)(q^{ab} - 1).
\eeq
\et

From Theorem~\ref{FP:theorem:polyIdentity}, 
we need only high school algebra to deduce the Sharp-Sylvester 
solution of the Frobenius problem.  
Recall that if $a_d \neq 0$ and 
\[
f(q) = a_d q^d + a_{d-1}q^{d-1} + \cdots + a_1 q + a_0
\]
is a polynomial of degree $d$, then the 
\emph{reciprocal polynomial}\index{reciprocal polynomial} of $f(t)$ is the polynomial 
\beq               \label{FP:reciprocal}
\hat{f}(q) = q^d f\left(\frac{1}{q}\right) 
= a_0 q^d + a_1q^{d-1} + \cdots + a_{d-1} q + a_d 
= \sum_{i=0}^d a_i q^{d-i}
\eeq
of degree at most $d$.  
For example, the degree 5 polynomial $f(q) = q^5 + q^2 + q$ has the degree 4 reciprocal polynomial 
$\hat{f}(q)  = q^4 + q^3 + 1$.

\bt                               \label{FP:theorem:Sylvester}
Let $A = \{a,b\}$, where $a$ and $b$ are distinct, relatively prime positive integers, 
and let $S(A) =  \{ai +bj: i,j \in \N_0\}$.  
\benum
\item[(i)] 
The Frobenius number of the set $A$ is  
\[
\mcf(A) = ab - a - b.
\]
\item[(ii)]
The numerical semigroup $S(A)$ is symmetric, and the genus of $A$  is 
\[
\mcg(A) = \frac{\mcf(A)+1}{2}  = \frac{(a-1)(b-1)}{2}.
\]
\eenum
\et

\begin{proof}
Because $\gcd(a,b) = 1$, at least one of the  integers $a$ and $b$ is odd, and so 
$(a-1)(b-1)/2$ is an integer.  

The degree of the polynomial $f_A(q)$ is the Frobenius number $\mcf(A)$, 
which is the largest integer not in $S(A)$.
Equating the degrees of the polynomials 
on the left and right sides of identity~\eqref{FP:polyIdentity} 
in Theorem~\ref{FP:theorem:polyIdentity}, we obtain 
\[
a+ b + 1 + \deg(f_A(q)) = 1 + ab
\]
and so
\[
\mcf(A) = \deg(f_A(q)) = ab - a -b.  
\]
It follows that the reciprocal polynomial of $f_A(q)$ is 
\[
\hat{f_A}(q) =  q^{ab-a-b}  f \left(  \frac{1}{q}  \right).
\]

Consider the polynomial 
\beq            \label{FP:fg}
g_A(q) = \sum_{ n =0}^{ab - a - b} q^n - f_A(q). 
\eeq

The polynomials on the left and right sides of identity~\eqref{FP:polyIdentity} 
have degree $ab + 1$.  
The reciprocal polynomial of the right side of~\eqref{FP:polyIdentity} is 
\[
q^{ab+1} \left(   \left(  \frac{1}{q} -1 \right)  \left(  \frac{1}{ q^{ab} } - 1 \right) \right)
= (q-1)(q^{ab} - 1).  
\]
The reciprocal polynomial of the left side of~\eqref{FP:polyIdentity} is 
\begin{align*}
q^{ab+1}&  \left(   \left(  \frac{1}{q^a} -1 \right)  \left(  \frac{1}{ q^b } - 1 \right)  
\left(   \left(  \frac{1}{q} -1 \right)  f \left(  \frac{1}{q}  \right) +1   \right)    \right)  \\
&  = (q^a - 1) (q^b - 1) \left( - (q-1)  q^{ab-a-b}  f \left(  \frac{1}{q}  \right) + q^{ab-a-b+1}    \right)  \\
&  = (q^a - 1) (q^b - 1) \left( - (q-1) \hat{f_A}(q)  + q^{ab-a-b+1}    \right).  
\end{align*}
Therefore, 
\beq                  \label{FP:polyIdentity-reciprocal}
(q^a - 1) (q^b - 1) \left( - (q-1) \hat{f_A}(q) +q^{ab-a-b+1}   \right)  = (q-1)(q^{ab} - 1).  
\eeq
Comparing identities~\eqref{FP:polyIdentity} and~\eqref{FP:polyIdentity-reciprocal}, 
we obtain 
\[
(q-1)f_A(q) + 1 =  - (q-1) \hat{f_A}(q) +q^{ab-a-b+1} 
\]
and so 
\[
(q-1) ( f_A(q) + \hat{f_A}(q)) =  q^{ab-a-b+1}  - 1 = (q-1) \left( 1 + q + q^2 + \cdots + q^{ab-a-b} \right) 
\] 
By identity~\eqref{FP:fg}, 
\[
 f_A(q) + \hat{f_A}(q) =  1 + q + q^2 + \cdots + q^{ab-a-b} = f_A(q) + g_A(q) 
\]
and so 
$
 \hat{f_A}(q) = g_A(q).  
$
Let 
\[
\varepsilon_n = \begin{cases}
1 & \text{ if $n \in S(A)$} \\
0 & \text{ if $n \notin S(A)$.} 
\end{cases}
\]
We have
\[
f_A(q) = \sum_{n =0}^{\mcf(A)} (1 -  \varepsilon_n) q^n
\]
and  
\[
g_A(q) = \sum_{n=0}^{\mcf(A)} \varepsilon_n q^n
\]
Recalling  formula~\eqref{FP:reciprocal} for the reciprocal polynomial, 
we obtain  
\[
 \hat{f_A}(q) = \sum_{n=0}^{\mcf(A)} (1 - \varepsilon_n ) q^{\mcf(A) - n} 
 = \sum_{n=0}^{\mcf(A)} \left(1 - \varepsilon_{\mcf(A) - n} \right)q^n = g_A(q).
\]
It follows that 
$
1 -  \varepsilon_{\mcf(A) - n} =   \varepsilon_n.
$
Equivalently, 
\[
 \varepsilon_n + \varepsilon_{\mcf(A) - n} = 1
\]
for all $n \in \{0,1,2,\ldots, \mcf(A) \}$.  
Therefore,  $n \in S(A)$ if and only if 
$ \varepsilon_n = 1$ if and only if $ \varepsilon_{\mcf(A) - n } = 0$ 
if and only if $\mcf(A) - n \notin S(A)$.  
Thus, the semigroup $S(A)$is symmetric, and the genus of $A$ is 
\[
\mcg(A)= f_A(1) = \frac{(a-1)(b-1)}{2}.  
\]
This completes the proof.  
\end{proof}

\section{A division algorithm in $E[x,y]$}
Let $E$ be a field, and let $E[t]$  and $E[x,y]$ be the polynomial rings 
in one and two variables, respectively.  
Let $A = \{ a,b\}$, where $a$ and $b$ are distinct, relatively prime positive integers, 
and let $S(A)  = \{ai +bj: i,j \in \N_0\}$. 
Consider the ring homomorphism $\Phi:E[x,y] \rightarrow E[t]$ defined by 
\beq              \label{FP:definePhi}
\Phi(x) = t^a \qqand \Phi(y) = t^b.   
\eeq 
For every polynomial $f(x,y) \in E[x,y]$, we have $\Phi(f(x,y) ) = f\left( t^a,t^b \right)$.  
Thus,  
\beq              \label{FP:definePhi-K}
\Phi\left( x^b - y^a \right) = \left( t^a \right)^b - \left(  t^b\right)^a = t^{ab} - t^{ab} = 0
\eeq
and so the kernel of $\Phi$ contains the polynomial $x^b - y^a$.  
We shall prove that the kernel of $\Phi$ is the principal ideal generated by $x^b - y^a$.

The image of $\Phi$ is the subring of $E[t]$ generated 
by the set $\left\{ t^a, t^b \right\}$, and denoted $E[ t^a, t^b]$.
Because $\Phi(x^iy^j) = t^{ai+bj}$, it follows that 
$\{ t^n: n \in S(A) \}$ is the set of monomials that appear in $E[ t^a, t^b]$, 
and so $E[ t^a, t^b]$ is a vector space over the field $E$ with basis $\{ t^n: n \in S(A) \}$.  

The proof of Theorem~\ref{FP:theorem:kernelPhi} 
uses the division algorithm for a  polynomial in $k$ variables 
by a sequence of $s$ polynomials in $k$ variables. 
(This is clearly explained in Cox-Little-O'Shea~\cite[Chapter 2]{cox-litt-oshe97}.)  
We need only the special case $k=2$ and $s=1$.
Lexicographically order the monomials $x^iy^j \in E[x,y]$ as follows:  
$x^{i_1} y^{j_1} \prec x^{i_2}y^{j_2}$ if $i_1 < i_2$ or if $ i_1 = i_2$ and $j_1 < j_2$.  
Every nonempty finite set of monomials has a largest element.  
The \emph{leading monomial} of a nonzero polynomial $f(x,y)$ is the largest monomial 
that appears in the polynomial with a nonzero coefficient.
For example, if $a$ and $b$ are positive integers, then 
$y^a = x^0 y^a \prec x^b y^0 = x^b$, and so $x^b$ is the leading monomial of 
the polynomial $x^b - y^a$ 
By the division algorithm, if the leading monomial of the polynomial $f(x,y)$ 
is $x^b$ for some positive integer $b$, then, for every polynomial $g(x,y) \in E[x,y]$,
there exists a ``quotient polynomial'' $q(x,y) \in E[x,y]$ and ``remainder polynomials'' 
$r_i(y) \in E[y]$ for $i=0,1,2,\ldots, b-1$ such that 
\beq           \label{FP:DivisionAlgorithm}
g(x,y) = q(x,y) f(x,y) +  \sum_{i=0}^{b-1} x^i  r_i(y).
\eeq

\bt               \label{FP:theorem:kernelPhi}
Let $a$ and $b$ be distinct, relatively prime positive integers. 
Let $E$ be a field, and consider the polynomial rings $E[x,y]$ and $E[t]$.
Define the ring homomorphism $\Phi: E[x,y] \rightarrow E[t]$ by~\eqref{FP:definePhi}.
The kernel of $\Phi$ is the principal ideal generated by $x^b - y^a$.  
\et

\begin{proof} 
Equation~\eqref{FP:definePhi-K} shows that $f (x,y) = x^b - y^a \in \kernel(\Phi)$.   

Let $g  \in \kernel (\Phi)$.  Using the division algorithm to divide $g(x,y)$ by $f(x,y)$, 
we obtain polynomials $q(x,y)$ and $r_0(y), r_1(y),\ldots, r_{b-1}(y)$ 
that satisfy~\eqref{FP:DivisionAlgorithm}.  
Let 
\[
r_i(y) = \sum_{j=0}^{d_i} c_{i,j} y^j
\]
for $i=0,1,2,\ldots, b-1$.  
Equation~\eqref{FP:DivisionAlgorithm} gives  
\begin{align*}
g(x,y) 
& = q(x,y) f(x,y) +  \sum_{i=0}^{b-1} \sum_{j=0}^{d_i} c_{i,j} x^i  y^j.
\end{align*}
Because $\Phi ( g(x,y) ) = \Phi( f(x,y) )= 0$, we obtain 
\begin{align*}
0  & = \Phi( g(x,y) )  
 = \Phi \left( q(x,y) f(x,y) +  \sum_{i=0}^{b-1} \sum_{j=0}^{d_i} c_{i,j} x^i  y^j \right) \\ 
& = \Phi ( q(x,y) )\Phi (  f(x,y) )+ \sum_{i=0}^{b-1} \sum_{j=0}^{d_i} c_{i,j} \Phi \left( x^i  y^j \right) \\ 
& =  \sum_{i=0}^{b-1} \sum_{j=0}^{d_i}  c_{i,j} t^{ai+bj}.
\end{align*}
For 
\[
(i_1,j_1) \in \{0,1,\ldots, b-1\} \times \N_0
\]
and
\[
(i_2,j_2) \in \{0,1,\ldots, b-1\} \times \N_0
\]
we have  
\[
t^{ai_1 + bj_1} = t^{a i_2 + b j_2 }
\]
if and only if 
\[
ai_1 + bj_1 = a i_2 + b j_2 
\]
if and only if 
\[
a(i_2 - i_1) = b(j_2-j_1).  
\]
Because $b$ divides $a(i_2 - i_1)$ and $\gcd(a,b) = 1$, it follows that 
$b$ divides $i_2 - i_1$.  
The inequality $0 \leq | i_2 - i_1 |\leq b-1$ implies  that 
$i_1 = i_2$ and so $j_1 = j_2$.   Thus, the integers in the set 
\[
\{ ai+bj: i \in \{0,1,\ldots, b-1\} \text{ and } j \in N_0 \} 
\]
are pairwise distinct elements of the semigroup $S(A)$, and  
the corresponding monomials $t^{ai+bj}$ are pairwise distinct.  
The polynomial identity 
\[
 \sum_{i=0}^{b-1} \sum_{j=0}^{d_i}  c_{i,j} t^{ai+bj} = 0
\]
implies that $c_{i,j} = 0$ for all $i$ and $j$, and so $r_i(y) = 0$
for $i = 0,1,\ldots, b-1$, and 
\[
g(x,y)  = q(x,y) f(x,y) .
\]
Thus, the kernel of $\Phi$ is the principal ideal generated by $f(x,y) $.   
This completes the proof.
\end{proof}

\section{Graded rings and modules}\label{FP:section:proof}

Let $R$ be a commutative ring with 1.  
The ring $R$ is \emph{graded}\index{graded ring}\index{ring!graded} 
if it contains a sequence $(R_n)_{n=0}^{\infty}$ of additive subgroups 
such that, first, as an additive abelian group, 
\[
R = \bigoplus_{n=0}^{\infty} R_n
\]
and, second, as a ring, multiplication in $R$ satisfies 
\beq                   \label{FP:RmRn}
R_m R_n \subseteq R_{m+n}
\eeq
for all $m,n \in \N_0$.  
In particular, $R_0 R_0 \subseteq R_0$ and $1 \in R_0$, and so $R_0$ is a ring.  
Similarly, for every $n \in \N_0$, we have $R_0R_n \subseteq R_n$, 
and so $R_n$ is an $R_0$-module.

Let $R = \bigoplus_{n=0}^{\infty} R_n$ be a graded ring. 
An $R$-module $M$ is  \emph{graded} 
if  $M$ contains a sequence $(M_n)_{n=0}^{\infty}$ of additive subgroups 
such that, first, as an additive abelian group, 
\[
M = \bigoplus_{n=0}^{\infty} M_n
\]
and, second, as an $R$-module, multiplication  satisfies 
\[
R_m M_n \subseteq M_{m+n}
\]
for all $m,n \in \N_0$.  
Because $R_0 M_n \subseteq M_n$, it follows that $M_n$ is an $R_0$-module 
for all $n \in \N_0$.  
If $f_n \in M_n$ for all $n$, and if $f = \sum_{n=0}^{\infty} f_n = 0$, 
then $f_n = 0$ for all $n$.  

If $R_0 = E$ is a field, then $M_n$ is a  vector space over $E$.  
If $M_n$ is a finite-dimensional vector space for all $n$, 
then the formal power series 
\[
H_M(q) = \sum_{n=0}^{\infty} \dim(M_n) q^n
\]
is called the \emph{Hilbert series} for $M$.

Relation~\eqref{FP:RmRn} implies that 
every graded ring $R$ is also a graded $R$-module with $M_n = R_n$ for all $n \in \N_0$.

Let $R = \bigoplus_{n=0}^{\infty} R_n$ be a graded ring, 
and let $M = \bigoplus_{n=0}^{\infty} M_n$ and $M' = \bigoplus_{n=0}^{\infty} M'_n$ 
be graded $R$-modules.  An $R$-module homomorphism $\Phi: M \rightarrow M'$ is 
\emph{graded}\index{graded $R$-module homomorphism} if 
$\Phi(M_n) \subseteq M'_n$ for all $n \in \N_0$.  
Define the $R_0$-module homomorphism $\varphi_n: M_n \rightarrow M'_n$ 
by restriction: $\varphi_n (f_n) = \Phi(f_n)$ for all $f_n \in M_n$.  
The kernel of $\varphi_n$ is a submodule, denoted $K_n$, of $M_n$, 
and so $K =  \bigoplus_{n=0}^{\infty} K_n$ is a graded $R$-module.  
If  $f = \sum_{n=0}^{\infty} f_n \in K$ with $f_n \in K_n$ for all $n$, then 
\[
\Phi(f) = \Phi\left( \sum_{n=0}^{\infty} f_n   \right) = \sum_{n=0}^{\infty} \Phi\left( f_n   \right) 
= \sum_{n=0}^{\infty} \varphi_n \left( f_n   \right) = 0
\]
and so $f \in \kernel(\Phi)$.  
Therefore, $K \subseteq \kernel(\Phi)$.  

Conversely, if $f \in \kernel(\Phi)$ and if 
$f = \sum_{n=0}^{\infty} f_n$ with $f_n \in M_n$ for all $n$, then 
\[
0 = \Phi(f) = \Phi\left( \sum_{n=0}^{\infty} f_n   \right) 
= \sum_{n=0}^{\infty} \varphi_n\left( f_n   \right) 
\]
with $ \varphi_n\left( f_n   \right) \in M_n$, and so $ \varphi_n\left( f_n   \right) = 0$.  
Therefore, $f_n \in \kernel(\varphi_n) = K_n$ and  $f \in K$.  
This proves that  $K = \kernel(\Phi)$.

Here are some examples of graded rings and modules.  
Let $E$ be a field.
The polynomial ring $E[t]$ is a vector space over $E$.  
For every $n \in \N_0$, let $R_n = Et^n$ be the one-dimensional subspace 
of $E[t]$ spanned by $t^n$.   
The identity $t^m t^n = t^{m+n}$ implies that $R_m R_n \subseteq R_{m+n}$, 
and so $E[t] = \bigoplus_{n=0}^{\infty} R_n$ is a graded ring with $\dim(R_n) = 1$ for all $n$.
As a graded $E[t]$-module, the Hilbert series for $E[t]$ is 
\[
H_{E[t]}(q) = \sum_{n=0}^{\infty} q^n = \frac{1}{1-q}.
\]

Let $a$ and $b$ be distinct, relatively prime positive integers, 
and let $S(A) = \{ai +bj: i,j \in \N_0\}$. 
Consider the ring $E[t^a,t^b]$.  
As a vector space over $E$, a basis for $E[t^a,t^b]$ is the set of monomials 
\[
\{t^{ai+bj} : (i,j) \in \N_0^2\} = \{t^n: n \in S(A)\}
\]
and so 
\[
E[t^a,t^b] =  \bigoplus_{n=0}^{\infty} R_n
\]
where
\beq   \label{FP:Eab}
R_n = \begin{cases}
E t^n & \text{if $n \in S(A)$} \\
0 & \text{ if $n \notin S(A)$.}  
\end{cases}
\eeq
Thus, $\dim(R_n) = 1$ if $n \in S(A)$ and $\dim(R_n) = 0$ if $n \notin S(A)$.   
Note that $R_0 = E$ because $0 \in S(A)$.  
As a graded $E[t^a, t^b]$-module, the Hilbert series for $E[t^a, t^b]$ is 
\[
H_{E[t^a, t^b]}(q) = \sum_{n \in S(A)} q^n 
= \sum_{n=0}^{\infty} q^n -  \sum_{n \in  \N_0 \setminus S(A)} q^n 
= \frac{1}{1-q} - f_A(q)
\]
where $f_A(q)$ is the polynomial defined by~\eqref{FP:f_A}.

A ring can be graded in many ways.  For example, in the polynomial ring $E[x,y]$, 
the degree of the monomial $x^iy^j$ is $i+j$.  
If $R_n$ is the vector subspace of $E[x,y]$ generated by the set of monomials 
of degree $n$, that is, by the set $\{x^n, x^{n-1}y, x^{n-2}y^2,\ldots, y^n\}$, 
then $R_n$ is an $E$-vector space of dimension $n+1$, 
and  $E[x,y] = \bigoplus_{n=0}^{\infty} R_n$ is a graded ring. 
With this grading by degree, as an $E[x,y]$-module, 
the Hilbert series for the polynomial ring $E[x,y]$ is 
\[
\sum_{n=0}^{\infty} (n+1) q^n = \frac{1}{(1-q)^2}.
\]

For the Frobenius problem, we use a different grading of $E[x,y]$.  
Let $E_n$ be the vector subspace of $E[x,y]$ generated by the set of monomials 
\[
\{x^i y^j : ai+bj = n\}.
\]
The number of monomials in this set is exactly $p_{a,b}(n)$, 
which is the number of partitions of $n$ 
into parts $a$ and $b$.  Euler observed that the generating function 
for this partition function is the formal power series 
\[
\sum_{n=0}^{\infty} p_{a,b}(n) q^n 
= \left( \sum_{i=0}^{\infty} q^{ai} \right) \left(  \sum_{j=0}^{\infty} q^{bj} \right) 
= \frac{1}{(1-q^a)(1-q^b)}.
\]
If $x^iy^j  \in E_m$ and $ x^k y^{\ell}  \in E_n$, then 
\begin{align*}
ai + bj & =  m \\ 
ak + b \ell & = n \\
a( i + k ) + b( j + \ell  ) & = m + n
\end{align*}
and so 
\[
 \left( x^iy^j  \right)  \left( x^k y^{\ell}  \right)= x^{ i + k} y^{ j +  \ell  }\in E_{m + n}.
\]
This implies that $E_m E_n \subseteq E_{m+n}$
and so $E[x,y] = \bigoplus_{n=0}^{\infty} E_n$ is a graded ring.  
Because $ai+bj = 0$ if and only if $i=j=0$,  
we have $E_0 = E$, and so $E_n$ is a vector space over the field $E$ with  
$\dim(E_n) = p_{a,b}(n)$.  
With this ``Frobenius grading,'', the Hilbert series for $E[x,y]$ is 
\beq            \label{Frobenius:grading}
H_{E[x,y]}(q) = \sum_{n=0}^{\infty} \dim(E_n) q^n 
= \sum_{n=0}^{\infty} p_{a,b}(n) q^n = \frac{1}{(1-q^a)(1-q^b)}.
\eeq

Let $E[t^a,t^b] =  \bigoplus_{n=0}^{\infty} R_n$,  
with $R_n$ defined by ~\eqref{FP:Eab}.    
We define a ``multiplication''  
\[
E[x,y] \times E[t^a,t^b]  \rightarrow E[t^a,t^b] 
\]
as follows: $x t^n = t^{a+n}$ and $y t^n = t^{b+n}$ for all $n \in S(A)$.  
Thus, if $x^iy^j \in E_m$ and $n \in S(A)$, then $ai+bj = m$
and $x^iy^j  t^n = t^{ai+bj+n} = t^{m+n}$. 
Therefore, $E_m R_n \subseteq R_{m+n}$ for all $m,n \in \N_0$, 
and $E[t^a,t^b]$ is a graded $E[x,y]$-module.

The function $\Phi:E[x,y] \rightarrow E[t^a,t^b]$ 
defined by $\Phi(x) = t^a$ and $\Phi(y) = t^b$
is a surjective ring homomorphism.   
Theorem~\ref{FP:theorem:kernelPhi} states that 
$K = \kernel(\Phi)$ is the principal ideal of $E[x,y]$ 
generated by the polynomial $f = f(x,y) = x^b - y^a \in E_{ab}$.   
Thus, the kernel of $\Phi$ is the graded ring  
\[
K = E[x,y] f = \bigoplus_{n=0}^{\infty} E_n f 
= \bigoplus_{n=ab}^{\infty} E_{n-ab} f = \bigoplus_{n=0}^{\infty} K_n
\]
where $E_{n-ab} f \subseteq E_n$ and 
\[
K_n = K \cap E_n 
= \begin{cases}
0 & \text{ if $n = 0,1,\ldots, ab-1$}\\
 E_{n-ab} f(x,y)  & \text{ if $n \geq ab$.}
\end{cases}
\]
For all $x^iy^j \in E_m$, we have $ai+bj = m$ and 
\[
 \Phi(x^iy^j) = t^{ai+bj}  = x^iy^j  \Phi(1)  \in E[t^a,t^b] 
\]
and so $\Phi:E[x,y] \rightarrow E[t^a,t^b]$ is also an $E[x,y]$-module homomorphism.  
The Hilbert series for  $K$ is  
\begin{align*}
H_K(q) 
& = \sum_{n=ab}^{\infty} \dim ( K_n ) q^n 
 = \sum_{n=ab}^{\infty} \dim (  E_{n-ab} f) q^n \\ 
& = \sum_{n=ab}^{\infty} \dim  (E_{n-ab} ) q^n 
=   q^{ab}   \sum_{n=ab}^{\infty} \dim (E_{n-ab} ) q^{n-ab} \\
& =   q^{ab}   \sum_{n=0}^{\infty} \dim (E_{n} ) q^{n} 
 = \frac{q^{ab}}{(1-q^a)(1-q^b)}. 
\end{align*}
The final equation comes from~\eqref{Frobenius:grading}.

We have the graded $E[x,y]$-modules $E[x,y] =  \bigoplus_{n=0}^{\infty} E_n$ 
and $E[t^a,t^b] =  \bigoplus_{n=0}^{\infty} R_n$, 
and the $E[x,y]$-module homomorphism $\Phi: E[x,y] \rightarrow E[t^a,t^b] $.  
The restriction of  $\Phi$ to $E_n$ is the linear transformation 
$\varphi:E_n \rightarrow R_n$, where $R_n \neq 0$ if and only if $n \in S(A)$.
If $n \in S(A)$, then there exist nonnegative integers $i$ and $j$ such that $n = ai+bj$.  
Because $x^i y^j \in E_n$ and $\varphi_n(x^i y^j) = t^{ai+bj } = t^n$, 
it follows that $\varphi_n$ is surjective for all $n \in \N_0$.  
The kernel of $\varphi_n$ is $K_n$.  
By the rank-nullity theorem in linear algebra, 
\beq          \label{FP:RankNullity}
\dim(E_n) = \dim(R_n) + \dim(K_n).
\eeq
Multiplying this equation by $q^n$ and summing over $n$, we obtain 
the following Hilbert series identity:  
\begin{align*}
H_{E[x,y]}(q) & = \sum_{n=0}^{\infty} \dim(E_n) q^n 
 = \sum_{n=0}^{\infty} ( \dim(R_n) + \dim(K_n) ) q^n \\
& = \sum_{n=0}^{\infty} \dim(R_n) q^n +   \sum_{n=0}^{\infty}  \dim(K_n) q^n \\
& = H_{E[t^a, t^b]} (q) + H_K(q).  
\end{align*}

Equivalently, \begin{align*}
 \frac{1}{(1-q^a)(1-q^b)}  
 = \frac{ 1}{1-q}  -  f_A(q) +  \frac{q^{ab}}{(1-q^a)(1-q^b)} = 0
\end{align*}
and so 
\[
(q^a - 1) (q^b - 1) (  (q-1) f_A(q)+1 ) = (q-1)(q^{ab} - 1).
\]
This completes the proof of Theorem~\ref{FP:theorem:polyIdentity}.

\def\cprime{$'$} \def\cprime{$'$} \def\cprime{$'$}
\providecommand{\bysame}{\leavevmode\hbox to3em{\hrulefill}\thinspace}
\providecommand{\MR}{\relax\ifhmode\unskip\space\fi MR }
\providecommand{\MRhref}[2]{%
  \href{http://www.ams.org/mathscinet-getitem?mr=#1}{#2}
}
\providecommand{\href}[2]{#2}

\end{document}